\documentclass[final,authoryear,1p,times]{elsarticle}

\RequirePackage[OT1]{fontenc}
\usepackage{amsmath,amsthm,amssymb}
\usepackage[isolatin]{inputenc}
\usepackage{amsfonts,dsfont}
\usepackage{pstricks}
\usepackage{comment}
\usepackage{natbib}
\usepackage{colordvi}
\usepackage{graphicx}

\newtheorem{theorem}{Theorem}[section]
\newtheorem{proposition}[theorem]{Proposition}

\newtheorem{lemma}[theorem]{Lemma}

\newcommand{\R}{\mathbb{R}} 

\renewcommand{\P}{\mathbb{P}}
\newcommand{\E}{\mathbb{E}} 
\DeclareMathOperator{\var}{Var} 
\DeclareMathOperator{\cov}{Cov} 
\newcommand{\telque}{\;|\;}
\newcommand{\FDR}{\mbox{FDR}}

\newcommand{\mbf}{\mathbf}

\newcommand{\ind}[1]{{\mbf{1}\{#1\}}}

\newcommand{\G}{{\widehat{\mathbb{G}}}_{m}}
\newcommand{\Gzero}{{\widehat{\mathbb{G}}}_{0,m}}
\newcommand{\Gone}{{\widehat{\mathbb{G}}}_{1,m}}
\newcommand{\Gzerotilde}{{\widetilde{\mathbb{G}}}_{0,m}}
\newcommand{\Gonetilde}{{\widetilde{\mathbb{G}}}_{1,m}}

\usepackage{graphics}
\newcommand{\FDP}{\mbox{FDP}} 
\newcommand{\T}{\mathcal{T}} 
\newcommand{\TBH}{\mathcal{T}^{{\tiny \mbox{BH} }}} 
\newcommand{\TBHpoint}{\dot{\mathcal{T}}^{{\tiny \mbox{BH} }}} 
\newcommand{\Z}{\mathbb{Z}} 
\newcommand{\W}{\mathbb{W}} 
\newcommand{\B}{\mathbb{B}} 
\newcommand{\qnorm}{ {\Phi}^{-1} }

\journal{Statistics and Probability Letters}

\begin{document}

\begin{frontmatter}



\title{On the false discovery proportion convergence under Gaussian equi-correlation}

\author[label1]{Delattre, S.}
\author[label2]{Roquain, E.}
\address[label1]{University of Paris Diderot, LPMA}
\address[label2]{UPMC University of Paris 6, LPMA}

\begin{abstract}
We study the convergence of the false discovery proportion (FDP) of the Benjamini-Hochberg procedure in the Gaussian equi-correlated model, when the correlation $\rho_m$ converges to zero as the hypothesis number $m$ grows to infinity.
By contrast with the standard convergence rate $m^{1/2}$ holding under independence, this study shows that the FDP converges to the false discovery rate (FDR) at rate $\{\min(m,1/\rho_m)\}^{1/2}$ in this equi-correlated model. 
\end{abstract}


\end{frontmatter}


\section{Introduction}
\label{}

When testing simultaneously a large number $m$ of null hypotheses, a popular global type I error, that can be traced back to \cite{See1968}, is the false discovery proportion (FDP), defined as the ratio of the number of erroneous rejections to the number of rejections. 
The average of this random variable, called the false discovery rate (FDR, introduced by \cite{BH1995}), has been studied by a considerable number of works, see for instance \cite{Sar2008} and \cite{BR2009} for a review. However, studying the FDR is not sufficient to catch the full behavior of the FDP, for instance a FDR control does not prevent the FDP from having large variations. Therefore, some other studies aim to directly control the upper-quantile of the FDP distribution, see e.g. \cite{GW2006,LR2005}, or to directly compute the distribution of the FDP, either non-asymptotically \cite{Chi2008, RV2010}, or asymptotically \cite{GW2004}.
Recently, \cite{Neu2008,Neu2009} computed the asymptotic distribution of the FDP \textit{actually} achieved by the Benjamini-Hochberg (BH) procedure (and some other adaptive procedures) under independence of the $p$-values. It is proved that the FDP converges to the FDR at the parametric rate $\sqrt{m}$. 
Furthermore, \cite{Far2006} showed that this convergence is unchanged under a specific short-range dependency between the $p$-values.

In this paper, we are interested in studying the convergence of the FDP of the BH procedure in the model where the test statistics have exchangeable Gaussian errors, with equi-correlation $\rho$ (allowing for instance long-range dependencies). This model has become quite standard in multiple testing (see e.g. \cite{BKY2006,FDR2007}), as it is a very simple instance of dependent $p$-value model.
From an intuitive point of view, the test statistics can be seen as independent test statistics plus a disturbance variable whose importance depends on the value of $\rho$. 
When $\rho \in (0,1)$ is fixed with $m$ (and in the ``ideal" setting where the $p$-values under the alternative are all equal to zero), 
 \cite{FDR2007} proved that the FDP of the BH procedure converges to a non-deterministic random variable that still depends on the disturbance variable. 
When $\rho_m\rightarrow 0$, we show here that this disturbance variable has no effect on the limit of the FDP anymore, which equals $\pi_0 \alpha$ (where $\pi_0$ is the proportion of true nulls), but
   can still have an effect on the asymptotic variance of the FDP or even on the convergence rate. More precisely, when $\rho_m\rightarrow 0$ our main result states
 that $\{\min(m,1/\rho_m)\}^{1/2} (\FDP_m - \pi_0\alpha) \leadsto \mathcal{N}(0,V)$ holds for a given $V>0$; in comparison with the independent case, we may distinguish the two following cases, recovering all the possible convergence regimes of $\rho_m$ to zero:
 \begin{itemize}
 \item[\textbullet] when $\lim_m m\rho_m =\theta \in [-1,+\infty)$, the limit of the FDP and the convergence rate   are the same as in the independent case. The asymptotic variance $V$ is larger if $\theta>0$, smaller if $\theta<0$ and is the same whenever $\theta=0$ (i.e. $m\rho_m \rightarrow 0$). 
 \item[\textbullet] when $\lim_m m\rho_m =+\infty$ and $\lim_m \rho_m = 0$ the convergence rate is $\rho_m^{-1/2}$ instead of ${m}^{1/2}$. 
 \end{itemize}
 On the one hand, this shows that the FDP of the BH procedure is still well concentrated around $\pi_0 \alpha$ under weak equi-correlation such that $\rho_m = O(1/m)$. On the other hand, this puts forward that  the concentration of the FDP of the BH procedure around the FDR may be arbitrarily slow when $\rho_m\rightarrow 0$, which is a striking result that has not been reported before to the best of our knowledge. 
Therefore, our recommendation is that  the  BH procedure can be used under Gaussian equi-correlation when $\rho_m=O(1/m)$ (including the case of a negative equi-correlation) but should be used carefully as soon as $m\rho_m\rightarrow \infty$, as the actual convergence rate of the FDP to the FDR might be much slower. 

The paper is organized as follows: Section~\ref{sec:model} presents the model, the notation and the main result. The latter is proved in Section~\ref{secproof}, including a generalization to any ``regular" thresholding procedure, 
recovering the so-called $\pi_0$-adaptive procedures studied in \cite{Neu2008}. 
Finally, some further points in connexion with our methodology are discussed in Section~\ref{sec:discuss}. 




%

\section{Setting and main result}\label{sec:model}

We observe $X_i = \tau_i+Y_i$, $1\leq i \leq m,$  where the parameter of interest is $(\tau_i)_i \in \{0,\mu\}^m$  (for a given $\mu>0$)  and the (unobservable) error vector
$(Y_1, ...,Y_m)$ is an exchangeable Gaussian vector with $\E Y_1 = 0$ and $\var Y_1 = 1$. We let $\rho_m=\cov(Y_1,Y_2)\in [-(m-1)^{-1},1]$. 
We consider the problem of the one-sided testing of the null ``$\tau_i=0$" against the alternative ``$\tau_i=\mu$", simultaneously for any $1\leq i \leq m$. To test each null, we define the $p$-value $p_i=\Phi(X_i)$, where $\Phi(z)=\P(Z\geq z)$ is the standard Gaussian upper-tail function. The c.d.f. of each $p$-value is denoted by $G_0(t)=t$ under the null and by $G_1(t)=\Phi(\Phi^{-1}(t)-\mu)$ under the alternative. 
The number of true nulls is denoted by $m_0(m)=|\{i\telque \tau_i=0\}| $ and is assumed to be of the form $\lfloor m \pi_0\rfloor $ for a  given proportion of true null $\pi_0 \in (0,1)$ independent of $m$. The ``mixture" c.d.f. of the $p$-values is denoted by $G(t)=\pi_0G_0(t)+(1-\pi_0) G_1(t)$. 
Next, we define the e.c.d.f.'s
$\Gzero(t)= (m_0(m))^{-1} \sum_{i=1}^{m} \ind{\tau_i=0} \ind{p_i\leq t}$, $\Gone(t)= (m-m_0(m))^{-1} \sum_{i=1}^{m} \ind{\tau_i>0} \ind{p_i\leq t}$ and $\G(t)=  m^{-1} \sum_{i=1}^{m} \ind{p_i\leq t}$. 

Given a pre-specified level $\alpha\in (0,1)$, the procedure of \cite{BH1995} can be defined as the procedure rejecting the nulls corresponding to $p_i\leq \TBH(\G)$ where the (data-driven) threshold $\TBH(\G)$ is $\max\{t\in[0,1]\telque \G(t)\geq t/\alpha\}$. 
Next, the false discovery proportion at a given threshold $t\in [0,1]$ is defined as the proportion of true nulls among the hypotheses having a $p$-value smaller than or equal to $t$: 
$$
\FDP_m(t)=\frac{|\{1\leq i \leq m \telque \tau_i=0, p_i\leq t \}|}{|\{1\leq i \leq m \telque p_i\leq t \}|\vee 1} = \frac{m_0(m)}{m} \frac{\Gzero(t)}{\G(t)\vee m^{-1}},
$$
where $|\cdot|$ denotes the cardinality function.

We now state our main result.
\begin{theorem}\label{main-th}
There is a unique point $t^\star\in (0,1)$ such that $G(t^\star)=t^\star/\alpha$ and we have
\begin{itemize}
\item[(i)] if $\lim_m m\rho_m =\theta \in [-1,+\infty)$, then
\begin{equation}
\sqrt{m} \big( \FDP_m(\TBH(\G)) - \pi_0 \alpha \big) \leadsto \mathcal{N} \left(0, \pi_0\alpha^2 \frac{1-t^\star}{ t^\star} + \theta \frac{\pi_0^2\alpha^2}{2\pi (t^{\star})^2} e^{-( \qnorm(t^\star))^{2}} \right)\,;
\label{conv-LSU-1}
\end{equation}
\item[(ii)] if $\lim_m m\rho_m =+\infty$ and $\lim_m \rho_m = 0$, then
\begin{equation}
\rho_m^{-1/2} \big( \FDP_m(\TBH(\G)) - \pi_0 \alpha \big) \leadsto \mathcal{N} \left(0,  \frac{\pi_0^2\alpha^2}{2\pi (t^{\star})^2} e^{-( \qnorm(t^\star))^{2}} \right)\,.
\label{conv-LSU-2}
\end{equation}
\end{itemize}
\end{theorem}

\section{Proof of Theorem~\ref{main-th}}\label{secproof}

\subsection{A more general result}\label{secgenresult}

In what follows, we denote the space of functions  from $[0,1]$ to $\R$ which are right-continuous and with left-hand limits  (Skorokhod's space)
 by $D(0,1)$ and the space of continuous functions from $[0,1]$ to $\R$ by $C(0,1)$.
The method for proving our result relies on the methodology let down by \cite{Neu2008} which consider the case of a general threshold function $\T : D(0,1) \rightarrow [0,1]$ which is Hadamard differentiable at $G$, tangentially to $C(0,1)$ (see \cite{Vaart1998} for a formal definition). The proof of Theorem~\ref{main-th} is presented here as a consequence of a more general theorem, true for any such threshold. The derivative of the threshold $\T$ at $G$, which is a continuous linear form on $C(0,1)$, is denoted by $\dot{\T}_G$. According to the Riesz representation theorem, the continuous linear form $\dot{\T}_G$ can be written as $\dot{\T}_G(F)=\int_0^1 F(t) \dot{\T}_G(dt)$, where we identified the linear form $\dot{\T}_G$  and the corresponding signed measure. 

\begin{theorem}\label{th-gen}
Let $\T : D(0,1) \rightarrow [0,1]$ be Hadamard differentiable at $G$, tangentially to $C(0,1)$, with derivative $\dot{\T}_G$. Let $q(t)=\pi_0 t/G(t)$ for $t>0$,  let $t^\star=\T(G)$ and assume $t^\star>0$. We set $\zeta_0 = \frac{q(t^\star)(1-q(t^\star)) }{t^\star} \delta_{t^\star} + \dot{q} (t^\star) \pi_0 \dot{\T}_G$, $\zeta_1 = -\frac{q(t^\star)(1-q(t^\star)) }{G_1(t^\star)} \delta_{t^\star} + \dot{q} (t^\star) (1-\pi_0) \dot{\T}_G$ and
 \begin{align*}
c(\T)=&\: (2\pi)^{-1/2}  \int_0^1 e^{-\frac{1}{2} ( \qnorm(t))^{2}} \zeta_0(dt)  +  (2\pi)^{-1/2}  \int_0^1 e^{-\frac{1}{2} ( \qnorm(t)-\mu)^{2}} \zeta_1(dt)\,;\\
\sigma^2(\T) =&\: \pi_0^{-1} \int_{[0,1]^2} (s\vee t -st) \zeta_0(ds)\zeta_0(dt) +(1-\pi_0)^{-1} \int_{[0,1]^2} (G_1(s\vee t) -G_1(s)G_1(t)) \zeta_1(ds)\zeta_1(dt)\,.
\end{align*}
Then the following holds:
\begin{itemize}
\item[(i)] if $\lim_m m\rho_m =\theta \in [-1,+\infty)$,
\begin{equation}
\sqrt{m} \big( \FDP_m(\T(\G)) - q(t^\star) \big) \leadsto \mathcal{N} \left(0, \sigma^2(\T)+ \theta c(\T)^2 \right)\,;
\label{conv-gen-1}
\end{equation}
\item[(ii)] if $\lim_m m\rho_m =+\infty$ and $\lim_m \rho_m = 0$,
\begin{equation}
\rho_m^{-1/2} \big( \FDP_m(\T(\G)) - q(t^\star) \big) \leadsto \mathcal{N} \left(0, c(\T)^2 \right)\,.
\label{conv-gen-2}
\end{equation}
\end{itemize}
 \end{theorem}

Let us now check that Theorem~\ref{main-th} follows from Theorem~\ref{th-gen}.
From \cite{Neu2008} Corollary 7.12, $\TBH:F\mapsto \max\{t\in[0,1]\telque F(t)\geq t/\alpha\}$ is  Hadamard differentiable at $G$, tangentially to $C(0,1)$, with derivative $\TBHpoint_G=(1/\alpha - \dot{G}(t^\star))^{-1} \delta_{t^\star}$. Moreover, $t^\star=\max\{t\in[0,1]\telque G(t)\geq t/\alpha\}$ is positive, because $\lim_{t\rightarrow 0^+} t/G(t)=0$.
Also, since $G(t^\star)=t^\star/\alpha$ and $ \dot{q} (t^\star) = (1/\alpha - \dot{G}(t^\star)) \pi_0 \alpha^2/ t^\star$, we may check that  $\zeta_1=0$ and $\zeta_0= (\pi_0\alpha/t^{\star}) \delta_{t^\star}$ in the above theorem, which leads to Theorem~\ref{main-th}.

\subsection{Proof of Theorem~\ref{th-gen}}

Let us now prove Theorem~\ref{th-gen}. First write $\FDP_m(\T(\G))=(m_0(m)/m) \pi_0^{-1} \Psi(\Gzero,\Gone)$, where for any $F_0,F_1$ in $D(0,1)$ with $F(\T(F))>0$ (letting $F=\pi_0 F_0+(1-\pi_0) F_1$), we put $\Psi(F_0,F_1) = \pi_0 \frac{F_0(\T(F))}{F(\T(F))}.$
From standard computations, $\Psi$ is Hadamard differentiable at $(G_0,G_1)$, tangentially to $C(0,1)^2$ and the derivative takes the form, for $(H_0,H_1) \in C(0,1)^2$, 
$\dot{\Psi}_{G_0,G_1}(H_0,H_1)= q(t^\star)(1-q(t^\star)) \left(\frac{H_0(t^\star)}{t^\star}- \frac{H_1(t^\star)}{G_1(t^\star)}\right) + \dot{q} (t^\star) \dot{\T}_G(H),$
where $H=\pi_0 H_0 + (1-\pi_0) H_1$. 
Applying the functional Delta method, this leads to the following useful result, which was essentially stated in \cite{Neu2008}.

\begin{proposition}\label{methodo-neu}
Let $\T : D(0,1) \rightarrow [0,1]$ be Hadamard differentiable at $G$, tangentially to $C(0,1)$, with derivative $\dot{\T}_G$. Let $q(t)=\pi_0 t/G(t)$ for $t>0$,  let $t^\star=\T(G)$ and assume $t^\star>0$. If for a given sequence $a_m \rightarrow \infty$ with $a_m=o(m)$,  
\begin{align}
a_m
 \left(\begin{array}{c}
\Gzero - G_0\\
\Gone - G_1
   \end{array} \right) 
   \leadsto 
     \left(\begin{array}{c}
   \W_0\\  \W_1   \end{array} \right) 
  \label{conv_process}
 ,
  \end{align}
where the convergence in distribution is relative to the Skorokhod topology and where $\W_0$ and $\W_1$ are processes with  continuous paths, then we have 
\begin{equation}
a_m \big( \FDP_m(\T(\G)) -q(t^\star) \big) \leadsto X
,\label{delta-method}
\end{equation}
where $X=\zeta_0(\W_0)+\zeta_1(\W_1)$ and $\zeta_0$, $\zeta_1$ are defined as in Theorem~\ref{th-gen}.
\end{proposition}

A convergence of the type \eqref{conv_process} in the particular Gaussian equi-correlated model is stated in Lemma~\ref{lem_conv-equi}.  
Using Proposition~\ref{methodo-neu}, this proves that \eqref{delta-method} holds both in the cases (i) and (ii) with $a_m=\sqrt{m}$ and $a_m=\rho_m^{-1/2}$, which respectively leads to \eqref{conv-gen-1} and \eqref{conv-gen-2} (the variance computations are straightforward). 

\subsection{Convergence of the e.c.d.f.'s in the Gaussian equi-correlated model}\label{convempprocess}

\begin{lemma}\label{lem_conv-equi}
Let $(\Z_0,\Z_1,Z)$ be a random variable such that $\Z_0\stackrel{(d)}{=} \pi_0^{-1/2} \B$, $\Z_1\stackrel{(d)}{=} (1-\pi_0)^{-1/2} \B\circ G_1$, $\B$ being a standard Brownian bridge on $[0,1]$, $\Z_0$ is independent from $\Z_1$, $Z\sim \mathcal{N}(0,1)$,  $\cov(Z,\Z_0(t))=  (2\pi)^{-1/2} \exp{(-\{\qnorm(t)\}^{2}/2)}$  and  $\cov(Z,\Z_1(t))  =  (2\pi)^{-1/2} \exp{(-\{\qnorm(t)-\mu\}^{2}/2)}$. Let also $U\sim\mathcal{N}(0,1)$ be independent of the vector $(\Z_0,\Z_1,Z)$.
Then we have the following convergences in law for the Skorokhod topology:
\begin{itemize}
\item[(i)] if $\lim_m m\rho_m =\theta \in [-1,+\infty)$, 
\begin{align}
\sqrt{m}
 \left(\begin{array}{c}
\Gzero - G_0\\
\Gone - G_1
   \end{array} \right) 
   \leadsto 
     \left(\begin{array}{c}
   \Z_0 +  (Z-\sqrt{1+\theta}U)\: \dot{\Phi}\circ\Phi^{-1}\\  \Z_1+(Z-\sqrt{1+\theta}U) \: \dot{\Phi}\circ(\Phi^{-1}-\mu)   \end{array} \right) \,;
  \label{conv_process_equi_cas1}
  \end{align}
  \item[(ii)] if $\lim_m m\rho_m =+\infty$ and $\lim_m \rho_m = 0$,
\begin{align}
\rho_m^{-1/2}
 \left(\begin{array}{c}
\Gzero - G_0\\
\Gone - G_1
   \end{array} \right) 
   \leadsto 
     \left(\begin{array}{c}
     U \: \dot{\Phi}\circ\Phi^{-1}\\   U  \: \dot{\Phi}\circ(\Phi^{-1}-\mu)  \end{array} \right) 
  \label{conv_process_equi_cas2}
 \,.
  \end{align}
  \end{itemize}
\end{lemma}

To prove Lemma~\ref{lem_conv-equi}, first remark that the distribution of the $X_i$'s may be realized as $X_i = \sqrt{1-\rho_m} (\xi_i - \overline{\xi}) + \sqrt{(1+(m-1)\rho_m)/m} \:U + \mu \ind{\tau_i>0},$ where $(\xi_1,...,\xi_m, U)$ are all i.i.d. $\mathcal{N}(0,1)$ variables and $\overline{\xi}$ denotes the empirical mean of the $\xi_i$'s. Let $\Gzero'(t)= (m_0(m))^{-1} \sum_{i : \tau_i=0}  \ind{\Phi(\xi_i)\leq t } $, $\Gone'(t)=(m-m_0(m))^{-1} \sum_{i:\tau_i>0} \ind{\Phi(\xi_i+\mu)\leq t }$
and
$$f_m(t,U,\rho_m) = (1-\rho_m)^{-1/2} \left(\Phi^{-1}(t)-  \sqrt{(1+(m-1)\rho_m)/m} \: U\right).$$
 The process $(\Gzero - G_0, \Gone - G_1)$ is then equal to  $V_m + W_m$ where 
 \begin{align*}
V_m(t)=  &     \left(\begin{array}{c}
(\Gzero'-G_0)(\Phi(f_m(t,U,\rho_m) + \overline{\xi})) \\
(\Gone'-G_1)(\Phi(f_m(t,U,\rho_m) + \overline{\xi}- \mu (1-\rho_m)^{-1/2}+\mu  ))    \end{array} \right) \\
 W_m(t)  =&
      \left(\begin{array}{c}
 \Phi(f_m(t,U,\rho_m) + \overline{\xi}) - t\\
 \Phi(f_m(t,U,\rho_m) + \overline{\xi}-\mu(1-\rho_m)^{-1/2})-  \Phi(\Phi^{-1}(t)-\mu)
   \end{array} \right) .
  \end{align*}
   
   Next, applying Donsker's theorem, we derive 
   $
   \sqrt{m} (\Gzero' - G_0,\Gone' - G_1,\overline{\xi})  {\leadsto} ( \Z_0,\Z_1 , Z)
   $,
 where $(\Z_0,\Z_1,Z)$ is defined as in Lemma~\ref{lem_conv-equi}.
  Since $\rho_m\rightarrow 0$, the inverse functions of $t\mapsto\Phi(f_m(t,U,\rho_m) + \overline{\xi})$ and $t\mapsto\Phi(f_m(t,U,\rho_m) + \overline{\xi}- \mu (1-\rho_m)^{-1/2}+\mu  )$ converge uniformly on $[0,1]$ to the identity a.s. Therefore, applying the Skorokhod's representation theorem, we get 
 \begin{align}
 \sqrt{m} (V_m,\overline{\xi}) \:{\leadsto}\: (\Z_0,\Z_1 ,Z ).
  \label{donsker}
  \end{align}
    
  Let us now consider the case (i), in which $\lim_m m\rho_m =\theta \in [-1,+\infty)$. In that case, a standard reasoning involving Taylor expansions of $\Phi$ and $y\mapsto\Phi(y \Phi^{-1}(t))$ 
     leads to 
    $$
    \sqrt{m} W_m(t) =     \left(\begin{array}{c}
 \dot{\Phi}(  \Phi^{-1}(t) ) (\sqrt{m}\: \overline{\xi}- \sqrt{1+\theta} \: U )  \\
 \dot{\Phi}(  \Phi^{-1}(t) -\mu ) (\sqrt{m} \:\overline{\xi}- \sqrt{1+\theta} \: U ) 
   \end{array} \right) 
   +
   \left(\begin{array}{c}
    R_{0,m}(t) \\
 R_{1,m}(t) 
   \end{array} \right),
    $$
    with remainder terms satisfying  $||R_{0,m}||_{\infty} \vee   ||R_{1,m}||_{\infty}  \rightarrow 0$ in probability. 
Since $U$ is independent of all the other variables, we derive from  \eqref{donsker} that 
  $
   \sqrt{m} (V_m, W_m) \:{\leadsto}\: \big(  \Z_0, \Z_1 ,  (Z-\sqrt{1+\theta}U)) \: \dot{\Phi}\circ\Phi^{-1} , (Z-\sqrt{1+\theta}U)) \:\dot{\Phi}\circ(\Phi^{-1}-\mu) \big)
.
$
This implies \eqref{conv_process_equi_cas1}. 
Consider now the case (ii), in which $\lim_m m\rho_m =+\infty$ and $\lim_m \rho_m = 0$. In that situation, we deduce from \eqref{donsker} that $\rho_m^{-1/2} V_m$ converges in probability to $0$. Furthermore, using that $\rho_m^{-1/2} \overline{\xi}$ tends to zero in probability, we obtain  that 
  $$\rho_m^{-1/2} W_m(t) =  \left(\begin{array}{c}
 \dot{\Phi}(  \Phi^{-1}(t) )  (-  U)  \\
 \dot{\Phi}(  \Phi^{-1}(t) -\mu ) (- U ) 
   \end{array} \right) 
   +
   \left(\begin{array}{c}
    T_{0,m}(t) \\
 T_{1,m}(t) 
   \end{array} \right),
   $$
with remainder terms satisfying  $||T_{0,m}||_{\infty} \vee   ||T_{1,m}||_{\infty}  \rightarrow 0$ in probability. This implies \eqref{conv_process_equi_cas2}. 


\section{Discussion: FDP convergence in the case $\rho_m=\rho\in (0,1)$}\label{sec:discuss}

\label{caseFDR2007}

When $\rho_m =\rho\in (0,1)$, we cannot expect that the FDP concentrates around the FDR as in Theorem~\ref{main-th} (see e.g. \cite{FDR2007} Theorem 2.1). 
As a consequence, even if the FDP has a mean below $\pi_0\alpha$ (because the false discovery rate of the BH procedure is below $\pi_0 \alpha$ for each $m$ for PRDS statistics, see Theorem~1.2 in \cite{BY2001}), the FDP can exceed $\pi_0\alpha + \varepsilon$ ($\varepsilon>0$) with a probability that does not vanish when $m$ grows to infinity.

We claim here that in the ideal situation where the parameters of the model $\pi_0$, $\mu$, $\rho$ are perfectly known, it is possible to modify the $p$-values so that the FDP convergence to the FDR keeps the parametric convergence rate $\sqrt{m}$. For this, we replace each test statistic $X_i$  by
$
\widetilde{X}_i= \sqrt{{m}/({(m-1)(1-\rho))}} (X_i-\overline{X}+(1-\pi_0)\mu),
$
so that $(\widetilde{X}_1,...,\widetilde{X}_m)$ is a Gaussian vector with variances equal to $1$, equi-correlation $\widetilde{\rho}_m=-(m-1)^{-1}$ and means $\E \widetilde{X}_i = \sqrt{{m}/({(m-1)(1-\rho))}}  \tau_i$. We build the corresponding $p$-values by letting $\widetilde{p}_i=\Phi(\widetilde{X}_i)$, which are uniform under the null and have the c.d.f. $\widetilde{G}_{1,m}(t)= \Phi(\Phi^{-1}(t)-\widetilde{\mu}_m)$ for $\widetilde{\mu}_m=(m/(m-1))^{1/2}\mu(1-\rho)^{-1/2}$ under the alternative. Although the latter depends (slightly) on $m$, we easily check that our methodology applies using $\widetilde{G}_1(t)= \Phi(\Phi^{-1}(t)-\widetilde{\mu})$ for $\widetilde{\mu}=\mu(1-\rho)^{-1/2}$ and that 
the following convergence holds:
$$
\sqrt{m} \big( \widetilde{\FDP}_m - \pi_0 \alpha \big) \leadsto \mathcal{N} \left(0\:,\: \pi_0\alpha^2 \frac{1-t_\rho^\star}{ t_\rho^\star} - \frac{\pi_0^2\alpha^2}{2\pi (t_\rho^{\star})^2} e^{-( \qnorm(t_\rho^\star))^{2}} \right),
$$
where $\widetilde{\FDP}_m$ denotes the FDP of the BH threshold $\TBH$ used with the $p$-values $\widetilde{p}_i$'s and where $t_\rho^{\star}\in (0,1)$ is the unique point $t\in(0,1)$ satisfying $\pi_0 t +(1-\pi_0) \Phi(\Phi^{-1}(t)-{\mu}(1-\rho)^{-1/2}) =t/\alpha$ (which depends on $\rho$). Of course, while this $p$-value modification greatly improves the concentration of the FDP, this approach is oracle because  $\pi_0$, $\mu$, $\rho$ are generally unknown. A correct estimation of the model parameters within such a procedure stays an open issue.

\bibliographystyle{model2-names}
\bibliography{biblio}







\end{document}